\documentclass[final]{siamltex}

\usepackage[parfill]{parskip}    
\usepackage{amsmath,amssymb,latexsym,enumerate,mathrsfs}
\usepackage{hyperref}
\usepackage{array}
\usepackage{times}

\usepackage{epstopdf}
\usepackage{subfig}
\usepackage{graphicx}

\usepackage{moreverb}
\usepackage{marvosym}
\usepackage{color,soul}
\definecolor{lightblue}{rgb}{.90,.95,1}
\sethlcolor{yellow}

\newtheorem{remark}{Remark}

\newtheorem{assumption}{\bf Assumption}
\def\qed{\hfill $\Box$}

\title{On the hierarchical optimal control of a chain of distributed systems
\thanks{Version -- August 10, 2015.}}

\author{Getachew K. Befekadu\thanks{Department of Mechanical and Aerospace Engineering, University of Florida - REEF, 1350 N. Poquito Rd, Shalimar, FL 32579, USA ({\tt gbefekadu@ufl.edu}).}
        \and Eduardo L. Pasiliao\thanks{Munitions Directorate, Air Force Research Laboratory, 101 West Eglin Blvd, Eglin AFB, FL 32542, USA ({\tt pasiliao@eglin.af.mil}).}}

\begin{document}
\maketitle

\renewcommand{\thefootnote}{\arabic{footnote}}

\begin{abstract}
In this paper, we consider a chain of distributed systems governed by a degenerate parabolic equation, which satisfies a weak H\"{o}rmander type condition, with a control distributed over an open subdomain. In particular, we consider two objectives that we would like to accomplish. The first one being of a controllability type that consists of guaranteeing  the terminal state to reach a target set starting from an initial condition; while the second one is keeping the state trajectory of the overall system close to a given reference trajectory on a finite, compact time intervals. We introduce the following framework. First, we partition the control subdomain into two disjoint open subdomains that are compatible with the strategy subspaces of the {\it leader} and that of the {\it follower}, respectively. Then, using the notion of Stackelberg's optimization (which is a hierarchical optimization framework), we provide a new result on the existence of optimal strategies for such an optimization problem -- where the {\it follower} (which corresponds to the second criterion) is required to respond optimally, in the sense of {\it best-response correspondence} to the strategy of the {\it leader} (that is associated to the controllability-type criterion) so as to achieve the overall objectives. Finally, we remark on the implication of our result in assessing the influence of the reachable target set on the optimal strategy of the {\it follower} in relation to the direction of {\it leader-follower} and {\it follower-leader} information flows.
\end{abstract}



\pagestyle{myheadings}
\thispagestyle{plain}
\markboth{G. K. BEFEKADU AND EDUARDO L. PASILIAO}{On the hierarchical optimal control}

\section{Introduction}  \label{S1}
\begin{figure}[bht]
\begin{center}
\includegraphics[width=110mm]{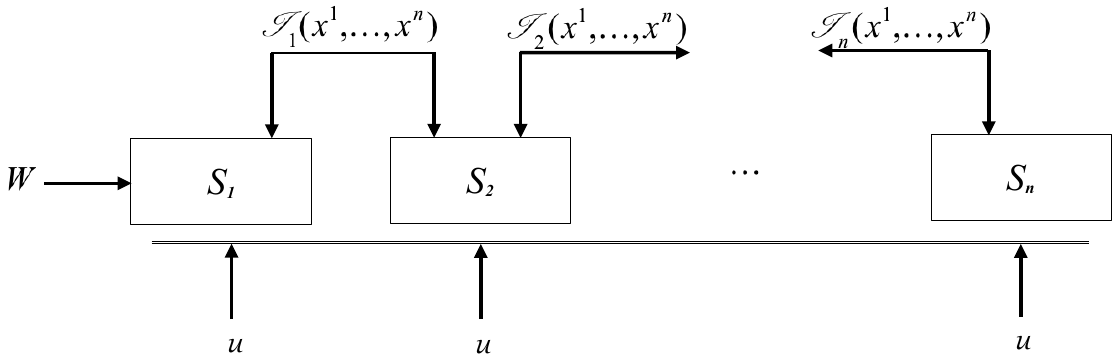}\\
{\small $ \begin{array}{c}\\
S_1:  \, \, d x_t^1 = f_1\bigl(t, x_t^1, \ldots,  x_t^n\bigr) dt + \sigma\bigl(t, x_t^1, \ldots, x_t^n\bigr)dW, \\
S_j: \,\, d x_t^j = f_j\bigl(t, x_t^{j-1}, \ldots,  x_t^n\bigr) dt, \quad j = 2, \ldots n,\\
     u \,\, \text{\it is a control distributed over an open subdomain}, \\
   \mathscr{I}_1, \mathscr{I}_2, \ldots, \mathscr{I}_n \,\, \text{\it are information for interconnecting subsystems}
\end{array}$}
\vspace{-2 mm}
\caption{A chain of distributed systems with random perturbations} \label{FigureDCS}
\vspace{-5 mm}
\end{center}
\end{figure}

In this paper, we consider the following distributed system, which is formed by a chain of $n$ subsystems (where $n \ge 2$), with a random perturbation that enters in the first subsystem and is then subsequently transmitted to other subsystems (see Figure~\ref{FigureDCS})
\begin{eqnarray}
\left.\begin{array}{l}
d x_t^1 = f_1\bigl(t, x_t^1, \ldots,  x_t^n\bigr) dt + \sigma\bigl(t, x_t^1, \ldots, x_t^n\bigr)dW(t) \\
d x_t^2 = f_2\bigl(t, x_t^1, \ldots,  x_t^n\bigr) dt  \\
d x_t^3 = f_3\bigl(t, x_t^2, \ldots,  x_t^n\bigr) dt  \\
 \quad\quad\quad~~ \vdots  \\
d x_t^n = f_n\bigl(t, x_t^{n-1}, x_t^n\bigr) dt, \quad 0 \le t \le T
\end{array}\right\},  \label{Eq1.1} 
\end{eqnarray}
where
\begin{itemize}
\item $x^i$ is an $\mathbb{R}^{d}$-valued state for the $i$th subsystem, with $i \in \{1, 2, \ldots, n\}$, 
\item the functions $f_1 \colon (0, \infty) \times \mathbb{R}^{nd} \rightarrow \mathbb{R}^d$ and $f_{j} \colon (0, \infty) \times \mathbb{R}^{(n-j+2)d} \rightarrow \mathbb{R}^d$, for $j = 2, \ldots, n$, are uniformly Lipschitz, with bounded first derivatives,
\item $\sigma \colon [0, \infty) \times \mathbb{R}^{nd} \rightarrow \mathbb{R}^{d \times m}$ is Lipschitz with the least eigenvalue of $\sigma\,\sigma^T$ uniformly bounded away from zero, i.e., 
\begin{align*}
 \sigma\bigl(t, x_t^1, \ldots, x_t^n\bigr)\,\sigma^T\bigl(t, x_t^1, \ldots, x_t^n\bigr) \ge \lambda I_d , \quad \forall (x_t^1, \ldots, x_t^n) \in \mathbb{R}^{nd}, \quad \forall t \ge 0,
\end{align*}
for some $\lambda > 0$,
\item $W$ (with $W(0) = 0$) is an $m$-dimensional standard Wiener process.
\end{itemize}

\begin{remark} \label{R1}
Note that such a distributed system has been discussed in various applications (e.g., see \cite{BarPV01}, \cite{BodL08} and \cite{Soi94} and the references therein). For example, when $n=2$, the equation in \eqref{Eq1.1} can be used to describe stochastic Hamiltonian systems (e.g., see \cite{BodL08} or \cite{Soi94} for additional discussions).
\end{remark}

Let us introduce the following notation that will be useful later. We use bold face letters to denote variables in $\mathbb{R}^{nd}$, for instance, $\mathbf{0}$ stands for a zero in $\mathbb{R}^{nd}$ (i.e., $\mathbf{0} \in \mathbb{R}^{nd}$) and, for any $t \ge 0$, the solution $\bigl(x_t^1, x_t^2, \ldots,  x_t^n\bigr)$ to \eqref{Eq1.1} is denoted by $\mathbf{x}_t$. Moreover, for $\bigl(t, (x^{j-1}, \ldots, x^n)\bigr) \in (0, \infty) \times \mathbb{R}^{(n-j+2)d}$, $j = 2, \ldots, n$, the function $x^j \mapsto f_j \bigl(t, x^{j-1}, \ldots, x^n\bigr)$ is continuously differentiable with respect to $x^j$ and its derivative denoted by $\bigl(t, x^{j-1}, \ldots, x^n\bigr) \mapsto D_{x^j} f_j \bigl(t, x^{j-1}, \ldots, x^n\bigr)$.

Then, we can rewrite the distributed system in \eqref{Eq1.1} as
\begin{align}
  d \mathbf{x}_t = F (t, \mathbf{x}_t ) dt + G \sigma (t, \mathbf{x}_t ) dW(t),  \label{Eq1.2}
\end{align}
where $F = \bigl [f_1, f_2, \ldots, f_n\bigr ]$ is an $\mathbb{R}^{nd}$-valued function and $G = \bigl[ I_d, 0, \ldots, 0 \bigr ]^T$ stands for an $(nd \times d)$ matrix that embeds $\mathbb{R}^d$ into $\mathbb{R}^{nd}$.

Let $\Omega$ be a regular bounded open domain in $\mathbb{R}^{nd}$, with smooth boundary $\Gamma$. For an open subdomain $U$ of $\Omega$, we consider the following distributed control system, governed by a partial differential equation (PDE) of parabolic type, with a control distributed over $U$, i.e.,
\begin{eqnarray}
\left.\begin{array}{l}
\dfrac{\partial y}{\partial t} + \mathcal{L}_{t,\mathbf{x}} y = u\chi_U \quad \text{in} \quad (0, T) \times \Omega  \\
 y (0, \mathbf{x} ) = 0 \quad \text{on} \quad \Omega  \\
 y (t, \mathbf{x} ) = 0 \quad \text{for} \quad (t, \mathbf{x}) \in \Sigma \triangleq (0, T) \times \Gamma
\end{array}\right\},   \label{Eq1.3}
\end{eqnarray}
where $u(t, \mathbf{x}) \in L^2((0,\,T) \times U)$ is a control function, $\chi_U$ is a characteristic function of the subdomain $U$ and $\mathcal{L}_{t,\mathbf{x}}$ is a  second-order operator given by\footnote{$\mathbf{x}^{j-1} \triangleq (x^{j-1}, \ldots, x^n)$ for $j = 2, \ldots n$.}
\begin{align}
 \mathcal{L}_{t,\mathbf{x}} = \dfrac{1}{2} \operatorname{tr} \bigl(a(t, \mathbf{x}) D_{x^1}^2\bigr) + f_1(t, \mathbf{x}) D_{x^1} + \sum\nolimits_{j=2}^n f_j(t, \mathbf{x}^{j-1}) D_{x^j},  \label{Eq1.4}
\end{align}
with $a(t, \mathbf{x})=\sigma(t, \mathbf{x}) \sigma^T(t, \mathbf{x})$.

\begin{remark} \label{R2}
Note that, in \eqref{Eq1.1}, the random perturbation enters in the first subsystem through the diffusive part and is then subsequently transmitted to other subsystems. As a result, such a distributed system is described by an $\mathbb{R}^{nd}$-valued diffusion process, which is degenerate in the sense that the second-order operator associated with it is a degenerate parabolic equation. Moreover, we also assume that the distributed system in \eqref{Eq1.1} satisfies a weak H\"{o}rmander condition (e.g., see \cite{Hor67} or \cite[Section~3]{Ell73} for additional discussions). 
\end{remark}

In what follows, we assume that the following statements hold true for the distributed system in \eqref{Eq1.1}.

\begin{assumption} \label{AS1} ~\\\vspace{-3mm}
\begin{enumerate} [(a)]
\item The functions $f_1(t, \mathbf{x})$ and $f_{j}(t, \mathbf{x}^{j-1})$ for $j = 2, \ldots, n$ are bounded $C^{\infty}((0, \infty) \times \mathbb{R}^{nd})$ and $C^{\infty}((0, \infty) \times \mathbb{R}^{(n-j+2)d})$-functions, respectively, with bounded first derivatives. Moreover, $\sigma(t, \mathbf{x})$ and $\sigma^{-1}(t, \mathbf{x})$ are bounded $C^{\infty} \bigl((0, \infty) \times \mathbb{R}^{nd}\bigr)$-functions, with bounded first derivatives.
\item The second-order operator $\mathcal{L}_{t,\mathbf{x}}$ in \eqref{Eq1.4} is hypoelliptic in $C^{\infty}((0, \infty) \times \mathbb{R}^{nd})$ (e.g., see \cite{Hor67} or \cite{Ell73}).
\end{enumerate}
\end{assumption}

\begin{remark}  \label{R3}
In general, the hypoellipticity assumption is related to a strong accessibility property of controllable nonlinear systems that are driven by white noise (e.g., see \cite{SusJu72} concerning the controllability of nonlinear systems, which is closely related to \cite{StrVa72} and \cite{IchKu74}; see also \cite[Section~3]{Ell73}). From Part~(b) of the above assumption, the Jacobian matrices $D_{x^1} f_1(t, \mathbf{x})$ and $D_{x^{j-1}} f_j(t, \mathbf{x}^{j-1})$ for $j=2, \ldots, n$ are assumed to be nondegenerate uniformly in time and space (i.e., they satisfy H\"{o}lder conditions both with respect to time and second variables). 
\end{remark}

Here it is worth mentioning that some studies on the controllability of systems that are governed by parabolic equations have been reported in literature (e.g., see \cite{Lio94} in the context of Stackelberg optimization; and \cite{AraFS15} and \cite{GuiMR13} in the context of Stackelberg-Nash controllability-type problem).\footnote{Recently, the authors in \cite{BefA15} and \cite{DelM10} have also provided some results, but in different contexts, pertaining to a chain of distributed systems with random perturbations.} Note that rationale behind our framework follows in some sense the settings of these papers. However, to our knowledge, the problem of optimal control for a chain of distributed system governed by degenerate parabolic equations has not been addressed in the context of hierarchical argument, and it is important because it provides a mathematical framework that shows how a hierarchical optimization framework can be systematically used to obtain optimal strategies for the {\it leader} and that of the {\it follower} (distributed over an open subdomain) for a chain of distributed system with random perturbations.\footnote{In this paper, our intent is to provide a theoretical framework, rather than considering a specific numerical problem or application.}

The remainder of this paper is organized as follows. In Section~\ref{S2}, using the remarks made above in Section~\ref{S1}, we state the optimal control problem for a chain of distributed system. Section~\ref{S3} presents our main results -- where we introduce a hierarchical optimization framework under which the {\it follower} is required to respond optimally, in the sense of {\it best-response correspondence} to the strategy of the {\it leader} (and vice-versa) so as to achieve the overall objectives. This section also contains results on the controllability-type problem for such a distributed system. For the sake of readability, all proofs are presented in Section~\ref{S4}. Finally, Section~\ref{S5} provides further remarks.

\section{Problem Formulation} \label{S2}
In this paper, we consider two objectives that we would like to accomplish. The first one being of a controllability type that consists of guaranteeing the terminal state to reach a target set from an initial condition; while the second one is keeping the state trajectory of the overall system close to a given reference trajectory on a finite, compact time intervals. Such a problem can be stated as follow:

{\bf Problem}: Find an optimal control strategy $u^{\ast}(t, \mathbf{x}) \in L^2((0,\,T) \times U)$ (which is distributed over $U$) such that
\begin{enumerate} [(i)]
\item {\it The first objective}: Suppose that we are given a target point $y^{t_{g}}$ in $L^2(\Omega)$.
\begin{enumerate} []
\item Then, we would like to have 
\begin{align}
 y(T; u^{\ast}) \in y^{t_{g}} + \alpha B, \quad \alpha > 0, \label{Eq2.1}
\end{align}
where $y(t; u^{\ast})$ denotes the function $\mathbf{x} \mapsto y(t, \mathbf{x}; u^{\ast})$, $B$ is a unit ball in $L^2(\Omega)$ and $\alpha$ is an arbitrary small positive number.\footnote{Note that the condition on the terminal state in \eqref{Eq2.1} is associated with a controllability-type problem with respect to an initial condition $y(0, \mathbf{x}) = 0$ on $\Omega$ (e.g., see \cite{Lio88} for additional discussions).}
\end{enumerate}
\item {\it The second objective}: Suppose that we are given a reference trajectory $y^{r_{f}}(t,\mathbf{x})$ in $L^2((0,\,T) \times \Omega)$.
\begin{enumerate} []
\item Then, we would like to have the state  trajectory $y(t,\mathbf{x}; u^{\ast})$ not too far from the reference $y^{r_{f}}(t,\mathbf{x})$ for all $t \in (0, T)$.
\end{enumerate}
\end{enumerate}

In order to make the above problem more precise, we specifically consider the following hierarchical cost functionals:
\begin{align}
 J_1(u) & = \dfrac{1}{2}{\int\int}_{(0,\,T) \times U}  u^2 d\mathbf{x} dt \notag \\
            & \quad \text{ s.t.} \quad  y(T; u) \in y^{t_{g}} + \alpha B, \quad \alpha > 0 \label{Eq2.2}
\end{align}
and
\begin{align}
 J_2(u) =  \dfrac{1}{2}{\int\int}_{(0,\,T) \times \Omega} & \bigl( y(t; u) - y^{r_{f}}(t,\mathbf{x}) \bigr)^2 d\mathbf{x} dt \notag \\
            & \quad \quad + \dfrac{\beta}{2}{\int\int}_{(0,\,T) \times U} u^2 d\mathbf{x} dt, \quad \beta > 0. \label{Eq2.3}
\end{align}

Note that, in general, finding such an optimal strategy $u^{\ast} \in L^2((0,\,T) \times U)$ that minimizes simultaneously the above cost functionals in \eqref{Eq2.2} and  \eqref{Eq2.3} is not an easy problem. However, in what follows, we introduce the notion of Stackelberg's optimization \cite{VonSta34} (which is a hierarchical optimization framework), where we specifically partition the control subdomain $U$ into two open subdomains $U_1$ and $U_2$ (with $U_1 \cap U_2 =  \varnothing$) that are compatible with the strategy subspaces of the {\it leader} and that of the {\it follower}, respectively. That is,
\begin{align}
 U = U_1 \cup U_2 \,\, \text{up to a set of measurable} \,\, U, \label{Eq2.4}
\end{align}
where the strategy for the {\it leader} (i.e., $u_1$) is from the subspace $L^2((0,\,T) \times U_1)$ and the strategy for the {\it follower} (i.e., $u_2$) is from the subspace $L^2((0,\,T) \times U_2)$.

Note that if $\chi_{U_i}$, for $i = 1, 2$, denotes the characteristic function for $U_i$ and $u_i$ is the restriction of the distributed control $u$ to $L^2((0,\,T) \times U_i)$. Then, the PDE in \eqref{Eq1.3} can be rewritten as  
\begin{eqnarray}
\left.\begin{array}{l}
\dfrac{\partial y}{\partial t} + \mathcal{L}_{t,\mathbf{x}} y = u_1\chi_{U_1} + u_2\chi_{U_2} \quad \text{in} \quad (0, T) \times \Omega  \\
 y(0, \mathbf{x}) = 0 \quad \text{on} \quad \Omega  \\
 y(t, \mathbf{x}) = 0 \quad \text{for} \quad (t, \mathbf{x}) \in \Sigma 
\end{array}\right\},   \label{Eq2.5}
\end{eqnarray}
where $y(t, \mathbf{x}; u) = y(t, \mathbf{x}; (u_1, u_2))$, with $(u_1, u_2) \in L^2((0,\,T) \times U_1) \times L^2((0,\,T) \times U_2)$.

Suppose that the strategy for the {\it leader} $u_1 \in L^2((0,\,T) \times U_1)$ is given. Then, the problem of finding an optmal strategy for the {\it follower}, i.e., $u_2^{\ast} \in L^2((0,\,T) \times U_2)$, which minimizes the cost functional $J_2$ is then reduced to finding an optimal solution for
\begin{align}
\inf_{u_2 \in L^2((0,\,T) \times U_2)}  J_2(u_1, u_2) \label{Eq2.6}
\end{align}
such that 
\begin{align}
 u_2^{\ast} = \mathcal{R}(u_1) \label{Eq2.7}
\end{align}
for some unique map $\mathcal{R} \colon L^2((0,\,T) \times U_1) \rightarrow L^2((0,\,T) \times U_2)$. Note that if we substitute $u_2^{\ast} = \mathcal{R}(u_1)$ into \eqref{Eq2.5}, then the solution $y(t, \mathbf{x}; (u_1, \mathcal{R}(u_1)))$ depends uniformly on $u_1 \in L^2((0,\,T) \times U_1)$. Moreover, the controllability-type problem in \eqref{Eq2.2} is then reduced to finding an optimal solution for
\begin{align}
& \inf_{u_1 \in L^2((0,\,T) \times U_1)}  J_1(u_1) \notag \\
& \quad \text{ s.t.} \quad  y(T; (u_1, \mathcal{R}(u_1))) \in y^{t_{g}} + \alpha B. \label{Eq2.8}
\end{align}
In the following section, we provide a hierarchical optimization framework for solving the above problems (i.e., the optimization problems in \eqref{Eq2.6}, together with \eqref{Eq2.7} and \eqref{Eq2.8}). Note that, for a given $u_1 \in L^2((0,\,T) \times U_1)$, the optimization problem in \eqref{Eq2.6} has a unique solution on $L^2((0,\,T) \times U_2)$ (cf. Proposition~\ref{P1}). Moreover, the optimization problem in \eqref{Eq2.8} makes sense if $y(T; (u_1, \mathcal{R}(u_1)))$ spans a dense subset of $L^2(\Omega)$, when $u_1$ spans the subspace $L^2((0,\,T) \times U_1)$ (cf. Propositions~\ref{P2} and \ref{P3}).

\section{Main Results} \label{S3}
In this section, we present our main results -- where we introduce a framework under which the {\it follower} is required to respond optimally, in the sense of {\it best-response correspondence} to the strategy of the {\it leader} (and vice-versa) so as to achieve the overall objectives. Moreover, such a framework allows us to provide a new result on the existence of optimal strategies for such optimization problems pertaining to a chain of distributed system with random perturbations.

\subsection{On the optimality distributed system for the follower} \label{S3.1}
Suppose that, for a given {\it leader} strategy $u_1 \in L^2((0,\,T) \times U_1)$, if $u_2^{\ast} \in L^2((0,\,T) \times U_2)$, i.e., the strategy for the {\it follower}, is an optimal solution to \eqref{Eq2.6} (cf. \eqref{Eq2.3}). Then, such a solution is characterized by the following optimality condition
\begin{align}
{\int\int}_{(0,\,T) \times \Omega} \bigl(y - y^{r_{f}}\bigr)\hat{y} d\mathbf{x} dt + \beta {\int\int}_{(0,\,T) \times U_2} u_2^{\ast} \hat{u}_2 d\mathbf{x} dt = 0, \notag\\
  \forall \hat{u}_2 \in L^2((0,\,T) \times U_2),  \label{Eq3.1}
\end{align}
where $y$ and $\hat{y}$ are, respectively, the solutions to the following PDEs
\begin{eqnarray}
\left.\begin{array}{l}
\dfrac{\partial y}{\partial t} + \mathcal{L}_{t,\mathbf{x}} y = u_1\chi_{U_1} + u_2^{\ast} \chi_{U_2} \quad \text{in} \quad (0, T) \times \Omega  \\
 y(0, \mathbf{x}) = 0 \quad \text{on} \quad \Omega  \\
 y(t, \mathbf{x}) = 0 \quad \text{for} \quad (t, \mathbf{x}) \in \Sigma
\end{array}\right\}  \label{Eq3.2}
\end{eqnarray}
and
\begin{eqnarray}
\left.\begin{array}{l}
\dfrac{\partial \hat{y}}{\partial t} + \mathcal{L}_{t,\mathbf{x}} \hat{y} = u_2^{\ast} \chi_{U_2} \quad \text{in} \quad (0, T) \times \Omega  \\
 \hat{y}(0, \mathbf{x}) = 0 \quad \text{on} \quad \Omega  \\
 \hat{y}(t, \mathbf{x}) = 0 \quad \text{for} \quad (t, \mathbf{x}) \in \Sigma
\end{array}\right\}.   \label{Eq3.3}
\end{eqnarray}

Furthermore, if we introduce an adjoint state $p$ as follow
\begin{eqnarray}
\left.\begin{array}{l}
-\dfrac{\partial p}{\partial t} + \mathcal{L}_{t,\mathbf{x}}^{\ast} p = y - y^{r_f} \quad \text{in} \quad (0, T) \times \Omega  \\
 p(T, \mathbf{x}) = 0 \quad \text{on} \quad \Omega  \\
 p(t, \mathbf{x}) = 0 \quad \text{for} \quad (t, \mathbf{x}) \in \Sigma
\end{array}\right\},   \label{Eq3.4}
\end{eqnarray}
where $\mathcal{L}_{t,\mathbf{x}}^{\ast}$ is the adjoint operator for $\mathcal{L}_{t,\mathbf{x}}$. Then, we have the following result which characterizes the map $\mathcal{R}$ in \eqref{Eq2.7} (i.e., the optimality distributed system for the {\it follower}).

\begin{proposition}\label{P1}
Let $u_1 \in L^2((0,\,T) \times U_1)$ be given. Suppose that the following PDE
\begin{eqnarray}
\left.\begin{array}{l}
\dfrac{\partial y}{\partial t} + \mathcal{L}_{t,\mathbf{x}} y = u_1\chi_{U_1} - \dfrac{1}{\beta} p \chi_{U_2}, \quad \text{in} \quad (0, T) \times \Omega  \\
-\dfrac{\partial p}{\partial t} + \mathcal{L}_{t,\mathbf{x}}^{\ast} p = y - y^{r_f} \quad \text{in} \quad (0, T) \times \Omega  \\
 y(0, \mathbf{x}) = 0 \quad \text{on} \quad \Omega  \\
 p(T, \mathbf{x}) = 0 \quad \text{on} \quad \Omega  \\
 y(t, \mathbf{x}) = p(t, \mathbf{x}) = 0 \quad \text{for} \quad (t, \mathbf{x}) \in \Sigma\\
\end{array}\right\},   \label{Eq3.5}
\end{eqnarray}
admits a unique solution pair $\bigl(y(u_1), p(u_1) \bigr)$ (which also depends uniformly on $u_1 \in L^2((0,\,T) \times U_1)$). Then, the optimality distributed system for the follower is given by
\begin{align}
\mathcal{R}(u_1) &= -\dfrac{1}{\beta} p(u_1) \chi_{U_2} \notag \\
                            &\equiv u_2^{\ast}. \label{Eq3.6}
\end{align}
\end{proposition}

\begin{remark} \label{R4}
The above proposition states that if the strategy of the {\it leader} $u_1 \in L^2((0,\,T) \times U_1)$ is given. Then, the strategy for the {\it follower} $u_2^{\ast} = \mathcal{R}(u_1)$, which is responsible for keeping the state trajectory $y(t,\mathbf{x}; (u_1, \mathcal{R}(u_1)))$ close to the given reference trajectory $y^{r_f}(t, \mathbf{x})$ on the time intervals $(0, T)$, is optimal in the sense of {\it best-response correspondence}. Later, in Proposition~\ref{P2}, we provide an additional optimality condition on the strategy of the {\it leader}, when such a correspondence is interpreted in the context of hierarchical optimization framework.   
\end{remark}

\subsection{On the optimality distributed system for the leader} \label{S3.2}
In this subsection, we provide an optimality condition on the strategy of the {\it leader} in \eqref{Eq2.2}, when the strategy for the {\it follower} satisfies the optimality condition of Proposition~\ref{P1}.

For a given $\xi \in L^2(\Omega)$, let $\varphi $ and $\vartheta$ be unique solutions to the following PDE
\begin{eqnarray}
\left.\begin{array}{l}
-\dfrac{\partial \varphi}{\partial t} + \mathcal{L}_{t,\mathbf{x}}^{\ast} \varphi = \vartheta \quad \text{in} \quad (0, T) \times \Omega  \\
\dfrac{\partial \vartheta}{\partial t} + \mathcal{L}_{t,\mathbf{x}} \vartheta = - \dfrac{1}{\beta} \varphi \chi_{U_2} \quad \text{in} \quad (0, T) \times \Omega  \\
 \vartheta(0, \mathbf{x}) = 0 \quad \text{on} \quad \Omega  \\
 \varphi(T, \mathbf{x}) = \xi \quad \text{on} \quad \Omega  \\
 \varphi(t, \mathbf{x}) =  \vartheta(t, \mathbf{x})  = 0 \quad \text{for} \quad (t, \mathbf{x}) \in \Sigma\\
\end{array}\right\}.   \label{Eq3.9}
\end{eqnarray}
Next, define the following linear decompositions  
\begin{align}
y = y_0 + z \quad \text{and} \quad p = p_0 + q \label{Eq3.10}
\end{align}
such that $y_0$ and $p_0$ are the unique solutions to the following PDE
\begin{eqnarray}
\left.\begin{array}{l}
\dfrac{\partial y_0}{\partial t} + \mathcal{L}_{t,\mathbf{x}} y_0 =  - \dfrac{1}{\beta} p_0 \chi_{U_2} \quad \text{in} \quad (0, T) \times \Omega  \\
-\dfrac{\partial p_0}{\partial t} + \mathcal{L}_{t,\mathbf{x}}^{\ast} p_0 = y_0 - y^{r_f} \quad \text{in} \quad (0, T) \times \Omega  \\
 y_0(0, \mathbf{x}) = 0 \quad \text{on} \quad \Omega  \\
 p_0(T, \mathbf{x}) = 0 \quad \text{on} \quad \Omega  \\
 y_0(t, \mathbf{x}) = p_0(t, \mathbf{x}) = 0 \quad \text{for} \quad (t, \mathbf{x}) \in \Sigma\\
\end{array}\right\}.   \label{Eq3.11}
\end{eqnarray}
Note that, from \eqref{Eq3.5} and \eqref{Eq3.11} together with \eqref{Eq3.10}, it is easy to show that $z$ and $q$ are the unique solutions to the following PDE 
\begin{eqnarray}
\left.\begin{array}{l}
\dfrac{\partial z}{\partial t} + \mathcal{L}_{t,\mathbf{x}} z = u_1^{\ast} \chi_{U_1} - \dfrac{1}{\beta} q \chi_{U_2} \quad \text{in} \quad (0, T) \times \Omega  \\
-\dfrac{\partial q}{\partial t} + \mathcal{L}_{t,\mathbf{x}}^{\ast} q = z \quad \text{in} \quad (0, T) \times \Omega  \\
 z(0, \mathbf{x}) = 0 \quad \text{on} \quad \Omega  \\
 q(T, \mathbf{x}) = 0 \quad \text{on} \quad \Omega  \\
 z(t, \mathbf{x}) = q(t, \mathbf{x})  = 0 \quad \text{for} \quad (t, \mathbf{x}) \in \Sigma\\
\end{array}\right\},   \label{Eq3.12}
\end{eqnarray}
where $u_1^{\ast} \in L^2((0,\,T) \times U_1)$ is an optimal strategy for the {\it leader} which satisfies additional conditions (see below \eqref{Eq3.14} and \eqref{Eq3.15}).

In what follows, let us denote the norm in $L^2(\Omega)$ by $\Vert \cdot \Vert_{L^2(\Omega)}$ and assume that $\xi \in L^2(\Omega)$ satisfies the following
\begin{align}
\bigl(z(T),\, \xi \bigr) = 0, \quad \forall u_1 \in L^2((0,\,T) \times U_1), \label{Eq3.13}
\end{align}
where $(\cdot,\, \cdot)$ denotes the scalar product in $L^2(\Omega)$. 

Then, we have the following result which characterizes the optimality condition for the {\it leader} in \eqref{Eq2.2}.
\begin{proposition} \label{P2}
The optimal strategy for the leader that minimizes 
\begin{align*}
& \inf_{u_1 \in L^2((0,\,T) \times U_1)}  J_1(u_1)\\
& \quad \text{ s.t.} \quad  y(T; (u_1, \mathcal{R}(u_1))) \in y^{t_{g}} + \alpha B 
\end{align*}
is given by
\begin{align}
u_1^{\ast} = \varphi (\xi) \chi_{U_1}, \label{Eq3.14}
\end{align}
where $\varphi(\xi)$ is given from the unique solution set $\bigl\{y(\xi), p(\xi), \varphi(\xi), \vartheta(\xi) \bigr\}$ for the optimality distributed system
\begin{eqnarray}
\left.\begin{array}{l}
\dfrac{\partial y}{\partial t} + \mathcal{L}_{t,\mathbf{x}} y = u_1^{\ast} \chi_{U_1} - \dfrac{1}{\beta} p \chi_{U_2} \quad \text{in} \quad (0, T) \times \Omega  \\
-\dfrac{\partial p}{\partial t} + \mathcal{L}_{t,\mathbf{x}}^{\ast} p = y - y^{r_f} \quad \text{in} \quad (0, T) \times \Omega  \\
-\dfrac{\partial \varphi}{\partial t} + \mathcal{L}_{t,\mathbf{x}}^{\ast} \varphi = \vartheta \quad \text{in} \quad (0, T) \times \Omega  \\
\dfrac{\partial \vartheta}{\partial t} + \mathcal{L}_{t,\mathbf{x}} \vartheta = - \dfrac{1}{\beta} \varphi \chi_{U_2} \quad \text{in} \quad (0, T) \times \Omega  \\
 y(0, \mathbf{x}) = \vartheta(0, \mathbf{x}) = 0 \quad \text{on} \quad \Omega  \\
 p(T, \mathbf{x}) = 0 \quad \text{on} \quad \Omega  \\
 \varphi(T, \mathbf{x}) = \xi \quad \text{on} \quad \Omega  \\
 y(t, \mathbf{x}) = p(t, \mathbf{x})=  \varphi(t, \mathbf{x}) =  \vartheta(t, \mathbf{x})  = 0 \quad \text{for} \quad (t, \mathbf{x}) \in \Sigma\\
\end{array}\right\}.   \label{Eq3.15}
\end{eqnarray}
Moreover, $\xi \in L^2(\Omega)$ is a unique solution to the following variational inequality\footnote{Note that, in \eqref{Eq3.16}, we write $y(T;\,\xi)$ to make explicitly the fact that the solution set $\bigl\{y(\xi), p(\xi), \varphi(\xi), \vartheta(\xi) \bigr\}$ of \eqref{Eq3.15} depends uniformly on $\xi \in L^2(\Omega)$.}
\begin{align}
\bigl(y(T; \xi) -y^{t_g},\, \hat{\xi} - \xi\bigr) + \alpha \bigl(\Vert \hat{\xi}\Vert_{L^2(\Omega)} - \Vert \xi \Vert_{L^2(\Omega)} \bigr) \ge 0, \quad \forall \hat{\xi} \in L^2(\Omega). \label{Eq3.16}
\end{align}
\end{proposition}

\begin{remark} \label{R5}
Note that the hierarchical optimization problem in Proposition~\ref{P2} requires the {\it follower} to respond optimally to the strategy of the {\it leader} in the sense of best-response correspondence, where such a correspondence is implicitly embedded in \eqref{Eq3.15} (see also Section~\ref{S4} for additional remarks).
\end{remark}

\subsection{On the controllability-type problem for the distributed system} \label{S3.3}
In the following, we consider the controllability-type problem in \eqref{Eq2.8}, where we provide a condition under which $y(T; (u_1, \mathcal{R}(u_1)))$ spans a dense subset of $L^2(\Omega)$, when $u_1$ spans the subspace $L^2((0,\,T) \times U_1)$.

\begin{proposition} \label{P3}
Suppose that Proposition~\ref{P2} holds true. Then, for every $y^{t_g} \in L^2(\Omega)$ and $\alpha > 0$ (which is arbitrary small), there exits $u_1 \in L^2((0,\,T) \times U_1)$ such that
\begin{align}
y(T; (u_1, \mathcal{R}(u_1))) \in y^{t_{g}} + \alpha B. \label{Eq3.30}
\end{align}
\end{proposition}
\begin{remark} \label{R6}
Note that the above proposition implicitly requires the strong accessibility property of the distributed system in \eqref{Eq1.1} which is concerned with the controllability property of nonlinear systems with random perturbations (see Assumption~\ref{AS1} and Remark~\ref{R3}).
\end{remark}

\begin{remark} \label{R7}
Following the same discussion as above (i.e., Subsections~\ref{S3.1}, \ref{S3.2} and \ref{S3.3}), we can also consider a family of hierarchical cost functionals $J_1(u_1), J_2(u_2), \dots, J_N(u_N)$, with $N \ge 3$, and a family of control strategies $\bigl\{u_i\bigr\}_{i =1}^N$ distributed over open subdomains $U_i \subset \Omega$, where $U = \cup_{i=1}^N U_i$, with $U_i \cap U_j = \varnothing$ for $i \neq j$. Moreover, if $u_N$ follows $u_{N-1}, \ldots, u_1$; and $u_{N-1}$ leads $u_N$ and, at the same time, it follows $u_{N-2}, \ldots, u_1$, etc. Then, for a given $u_1 \in L^2((0,\,T) \times U_1)$, such a whole hierarchical optimization problem could be solved if there exist some maps $\mathcal{R}_j \colon L^2((0,\,T) \times U_j) \rightarrow L^2((0,\,T) \times U_{j+1})$ (that also depend on the cost functionals) such that $u_{j+1}^{\ast} = \mathcal{R}_j(u_j)$ for $j=1, \ldots, N-1$ (e.g., see \cite[Section~3]{Lei78} for further discussions on one-leader and many-followers).
\end{remark}

\section{Proof of the Main Results} \label{S4}
In this section, we give the proofs of our results.

\subsection{Proof of Proposition~\ref{P1}} \label{S4.1}
For a given $u_1 \in L^2((0,\,T) \times U_1)$, let $y$ and $p$ be the unique solutions of \eqref{Eq3.5}. If we multiply the second equation in \eqref{Eq3.5} by $\hat{y}$ and integrate by parts. Further, noting the PDEs in \eqref{Eq3.3} and \eqref{Eq3.4}, then we have the following
\begin{align}
{\int\int}_{(0,T)\times \Omega} \bigl(y - y^{r_f}\bigr) \hat{y} d \mathbf{x}dt &= {\int\int}_{(0,T)\times \Omega} \biggl(-\dfrac{\partial p}{\partial t} + \mathcal{L}_{t,\mathbf{x}}^{\ast} p \biggr) \hat{y} d \mathbf{x}dt \notag \\
& = {\int\int}_{(0,T)\times \Omega} p \biggl(\dfrac{\partial \hat{y}}{\partial t} + \mathcal{L}_{t,\mathbf{x}} \hat{y} \biggr) d \mathbf{x}dt  \notag \\
& = {\int\int}_{(0,T)\times U_2} p u_2^{\ast} d \mathbf{x}dt.  \label{Eq3.7}
\end{align}
Moreover, using the optimality condition in \eqref{Eq3.1} together with \eqref{Eq3.7}, we obtain
\begin{align}
p \chi_{U_2} + \beta u_2^{\ast} = 0, \label{Eq3.8}
\end{align}
which further gives an optimal strategy for the {\it follower} as
\begin{align*}
 u_2^{\ast} &= -\dfrac{1}{\beta} p \chi_{U_2} \\
  &\equiv \mathcal{R}(u_1),
\end{align*}
where $p$ is from the unique solution set $\{p(u_1), y(u_1)\}$ of \eqref{Eq3.5} that depends uniformly on $u_1 \in L^2((0,\,T) \times U_1)$. This completes the proof of Proposition~\ref{P1}. \qed

\subsection{Proof of Proposition~\ref{P2}} \label{S4.2}
Note that the optimization problem for the {\it leader} in \eqref{Eq2.2} is equivalent to
\begin{align*}
  &\inf_{u_1} \dfrac{1}{2}{\int\int}_{(0,\,T) \times U_1}  u_1^2 d\mathbf{x} dt  \\
            & \quad \text{ s.t.} \quad  y(T; (u_1, \mathcal{R}(u_1))) \in y^{t_{g}} - y_0(T) + \alpha B \quad (\text{see}~\eqref{Eq3.10}).
\end{align*}
Introduce the following cost functionals
\begin{align}
\bar{J}_1(u_1)  = \dfrac{1}{2}{\int\int}_{(0,\,T) \times U_1}  u_1^2 d\mathbf{x} dt \label{Eq3.17}
\end{align}
and
\begin{align}
\bar{J}_2(u_1) = \biggl\{\begin{array}{l}
0 \quad\quad  \text{if} \quad \xi \in y^{t_{g}} - y_0(T) + \alpha B  \\
+\infty \quad \text{otherwise on} \,\, L^2(\Omega)
\end{array} \label{Eq3.18}
\end{align}
Let $\mathcal{H} \in \mathscr{L}(L^2((0,\,T) \times U_1); L^2(\Omega))$ be a bounded linear operator such that\footnote{$\mathscr{L}(L^2((0,\,T) \times U_1); L^2(\Omega))$ denotes a family of bounded linear operators.}
\begin{align}
\mathcal{H} u_1 = z(T; u_1). \label{Eq3.19}
\end{align}
Then, the optimization problem in \eqref{Eq2.2} is equivalent to
\begin{align}
  \inf_{u_1 \in L^2((0,\,T) \times U_1)} \biggl\{ \bar{J}_1(u_1)  + \bar{J}_2(u_1) \biggr\}. \label{Eq3.20}
\end{align}
Furthermore, using Fenchel duality theorem (e.g., see \cite{Roc67} or \cite{EkeT76}), we have the following
\begin{align}
  \inf_{u_1 \in L^2((0,\,T) \times U_1)} \biggl\{ \bar{J}_1(u_1)  + \bar{J}_2(u_1) \biggr\} = - \inf_{\xi \in L^2(\Omega)} \biggl\{ \bar{J}_1^{\ast}(\mathcal{H}^{\ast} \xi)  + \bar{J}_2^{\ast}(-\xi) \biggr\}, \label{Eq3.21}
\end{align}
where $\mathcal{H}^{\ast}$ is the adjoint operator of $\mathcal{H}$ and the conjugate functions $\bar{J}_i^{\ast}$ are  given by
\begin{align}
  \bar{J}_i^{\ast} (\varphi) = \sup_{\hat{\varphi}} \bigl\{ (\varphi, \hat{\varphi}) - \bar{J}_i(\hat{\varphi}) \bigr\}, \quad i = 1, 2. \label{Eq3.22}
\end{align}
Note that if we multiply the first equation (respectively, the second one) in \eqref{Eq3.15} by $z$ (respectively, by $q$) and integrate by parts, then we obtain the following
\begin{align}
 (z(T), \xi) = {\int\int}_{(0,\,T) \times U_1} \varphi u_1^{\ast} d\mathbf{x} dt.  \label{Eq3.23}
\end{align}
Then, for $\xi \in L^2(\Omega)$ that satisfies \eqref{Eq3.13}, we have the following 
\begin{align}
  \mathcal{H}^{\ast} \xi = \varphi \chi_{U_1}, \label{Eq3.24}
\end{align}
where $\varphi$ is from the unique solutions of \eqref{Eq3.15}. 

Note that
\begin{align}
  \bar{J}_1^{\ast} (u_1^{\ast}) = J_1(u_1^{\ast}) \label{Eq3.25}
\end{align}
and
\begin{align}
  \bar{J}_2^{\ast} (\xi) = (\xi, y^{t_{g}} - y_0(T)) + \alpha \Vert \xi \Vert_{L^2(\Omega)}. \label{Eq3.26}
\end{align}
Then, the optimization problem in \eqref{Eq2.2} is equivalent to 
\begin{align}
 & \inf_{\xi} \dfrac{1}{2}{\int\int}_{(0,\,T) \times U_1}  \varphi^2 d\mathbf{x} dt + \alpha \Vert \xi \Vert_{L^2(\Omega)} - (\xi, y^{t_{g}} - y_0(T)) \notag \\
  & \quad \text{ s.t.} \quad  y(T; (u_1, \mathcal{R}(u_1))) \in y^{t_{g}} - y_0(T) + \alpha B. \label{Eq3.27}
\end{align}
Let $\xi \in L^2(\Omega)$ be a unique solution to the following variational inequality
\begin{align}
{\int\int}_{(0,\,T) \times U_1} \varphi(\hat{\varphi} - \varphi) d\mathbf{x} dt  + \bigl(y(T; \xi) -y^{t_g},\, \hat{\xi} - \xi\bigr) +& \alpha \bigl(\Vert \hat{\xi}\Vert_{L^2(\Omega)} - \Vert \xi \Vert_{L^2(\Omega)} \bigr) \ge 0, \notag \\
&\quad \quad  \forall \hat{\xi} \in L^2(\Omega). \label{Eq3.28}
\end{align}
Moreover, if we multiply the first equation (respectively, the second one) in \eqref{Eq3.9} by $(\hat{\varphi} - \varphi)$ (respectively, by $(\hat{\vartheta} - \vartheta))$ and integrate by parts, we obtain the following
\begin{align}
 {\int\int}_{(0,\,T) \times U_1} \varphi(\hat{\varphi} - \varphi) d\mathbf{x} dt = \bigl(z(T), \hat{\xi} - f\bigr). \label{Eq3.29}
\end{align}
Thus, if we substitute \eqref{Eq3.29} into \eqref{Eq3.28}, then we obtain \eqref{Eq3.16}. This completes the proof of Proposition~\ref{P2}. \qed

\subsection{Proof of Proposition~\ref{P3}} \label{S4.3}
From Proposition~\ref{P2}, suppose that the PDE in \eqref{Eq3.15} admits unique solutions (i.e., $y(t,\mathbf{x}; \xi)$, $p(t,\mathbf{x}; \xi)$, $\varphi(t,\mathbf{x}; \xi)$ and $\vartheta(t,\mathbf{x}; \xi)$ for $(t, \mathbf{x})\in (0,\,T) \times \Omega$ that depend uniformly on $\xi \in L^2(\Omega)$). Then, noting \eqref{Eq3.13}, the condition in \eqref{Eq3.23} becomes
\begin{align}
 \varphi(t,\mathbf{x}) = 0 \quad \text{for} \quad (t, \mathbf{x}) \in (0,\,T) \times U_1,  \label{Eq3.31}
\end{align}
which implies the following conditions (see also \eqref{Eq3.9})
\begin{align*}
 \vartheta \chi_{U_1} &= 0, \\
\dfrac{\partial \vartheta}{\partial t} + \mathcal{L}_{t,\mathbf{x}} \vartheta &= 0 \quad \text{in} \quad (0,\,T) \times \bigl(U \setminus U_2\bigr)
\end{align*}
and
\begin{align*}
-\dfrac{\partial \varphi}{\partial t} + \mathcal{L}_{t,\mathbf{x}}^{\ast} \varphi = 0 \quad \text{in} \quad (0,\,T) \times \bigl(U \setminus U_2\bigr).
\end{align*}
Furthermore, using Mizohata's uniqueness theorem (e.g., see \cite{SauSc87}) together with the regularity conditions on $f_i(t, \mathbf{x}^{i-1})$ and $a(t, \mathbf{x})$ (see also Assumption~\ref{AS1}), then we obtain the following
\begin{align}
 \vartheta (t,\mathbf{x}) = 0  \quad \text{for} \quad (t, \mathbf{x}) \in (0,\,T) \times \bigl(U \setminus U_2\bigr), \label{Eq3.32}
\end{align}
which requires $\xi$ to have a zero value outside of $U_2$ (cf. \eqref{Eq3.24} and \eqref{Eq3.13}). 

Next, consider the restriction of $\varphi$ and $\vartheta$ to $(0,\,T) \times U_2$ such that
\begin{eqnarray}
\left.\begin{array}{l}
-\dfrac{\partial \varphi}{\partial t} + \mathcal{L}_{t,\mathbf{x}}^{\ast} \varphi = \vartheta \quad \text{in} \quad (0, T) \times U_2  \\
\dfrac{\partial \vartheta}{\partial t} + \mathcal{L}_{t,\mathbf{x}} \vartheta = - \dfrac{1}{\beta} \varphi \chi_{U_2} \quad \text{in} \quad (0, T) \times U_2  \\
 \varphi, \,\,  \dfrac{\partial \varphi}{\partial \mathbf{x}^i}, \,\, \vartheta, \,\, \dfrac{\partial \vartheta}{\partial \mathbf{x}^i} = 0, \quad \text{for} \quad (t, \mathbf{x}) \in (0, T) \times \partial \, U_2 \\
 \varphi(T, \mathbf{x}) = \xi \chi_{U_1}, \quad \vartheta(0, \mathbf{x}) = 0  \quad \text{on} \quad U_2 \\
\end{array}\right\}. \label{Eq3.33}
\end{eqnarray}
Then, from \eqref{Eq3.33}, it remains to show that $\xi \chi_{U_2} = 0$, which is a sufficient condition for $y(T; (u_1, \mathcal{R}(u_1)))$ to span a dense subset of $L^2(\Omega)$, when $u_1$ spans the subspace $L^2((0,\,T) \times U_1)$. 

Note that the first two equations in \eqref{Eq3.33} imply the following
\begin{align}
\left(\dfrac{\partial}{\partial t} + \mathcal{L}_{t,\mathbf{x}} \right) \left(-\dfrac{\partial}{\partial t} + \mathcal{L}_{t,\mathbf{x}}^{\ast} \right) \varphi + \dfrac{1}{\beta} \varphi = 0, \quad \text{on} \quad (0, T) \times U_2. \label{Eq3.34}
\end{align}
which is a quasi-elliptic equation; and in view of Cauchy problems on bounded domains (e.g., see \cite[Theorem~6.6.1]{StrVa79}), for any fixed $t \in (0, T)$, $\varphi(t, \mathbf{x})$ is analytic in $U_2$, with Cauchy data zero on $\partial \,U_2$. As a result of this, $\varphi(t, \mathbf{x}) = 0$ on $(0, T) \times \partial \,U_2$ and also continuous in $t$, then we have $\varphi(0, \mathbf{x}) = 0$ and $\varphi(t, \mathbf{x})=\vartheta(t, \mathbf{x})=0$ for $(t, \mathbf{x}) \in (0, T) \times \partial \, U_2$, which implies  $\xi \chi_{U_2} = 0$ (cf. \eqref{Eq3.33}, since $\varphi(T, \mathbf{x}) = \xi \chi_{U_1}$). This completes the proof of Proposition~\ref{P3}.  \qed

\section{Further remarks} \label{S5}
In this section, we briefly comment on the implication of our result in assessing the influence of the reachable target set on the strategy of the {\it follower} in relation to the direction of {\it leader-follower} and {\it follower-leader} information flows.

Note that the statement in Proposition~\ref{P1} (i.e., the optimality distributed system for the {\it follower}) is implicitly accounted in Proposition~\ref{P2} (cf. \eqref{Eq3.15}). Hence, the optimal strategy for the {\it follower} (cf. \eqref{Eq3.6}) is given by  
\begin{align*}
u_2^{\ast} &= - \dfrac{1}{\beta} p(\xi) \chi_{U_2}  \\
       &\equiv \mathcal{R}(u_1),
\end{align*}
where $\xi \in L^2(\Omega)$ is a minimum solution to the variational inequality in \eqref{Eq3.16} and it also assumes a zero value outside of $U_2$ (i.e., $\xi \chi_{U_2} = 0$). Moreover, such a minimum solution lies in a certain dense subset of $L^2(\Omega)$, which is spanned by $y(T; (u_1, \mathcal{R}(u_1)))$, when $u_1$ spans $L^2((0,\,T) \times U_1)$ (cf. Proposition~\ref{P3}).

Note that, from Proposition~\ref{P2} (cf. \eqref{Eq3.24}), we also observe that the optimal strategy for the {\it leader} is given by
\begin{align*}
u_1^{\ast} = \varphi(\xi) \chi_{U_1},
\end{align*}
which is implicitly conditioned by the target set $y^{t_{g}} + \alpha B$, where $\alpha$ is an arbitrary small positive number (cf. \eqref{Eq2.2} or \eqref{Eq2.8}). Moreover, the terminal state is guaranteed to reach the target set starting from an initial condition $y(0, \mathbf{x}) = 0$ on $\Omega$ (i.e., $y(T; (u_1^{\ast}, \mathcal{R}(u_1^{\ast}))) \in y^{t_{g}} + \alpha B$); and the state trajectory $y(t,\mathbf{x}; (u_1^{\ast}, \mathcal{R}(u_1^{\ast})))$ is not too far from the reference $y^{r_{f}}(t,\mathbf{x})$ for all $t \in (0, T)$. As a result of this, such interactions constitute a constrained information flow between the {\it leader} and that of the {\it follower} (i.e., an information flow from {\it leader-to-follower} and vice-versa) that captures implicitly the influence of the reachable target set on the strategy of the {\it follower}.

\end{document}